\documentclass[12pt,letterpaper,oneside,reqno]{amsart}

\usepackage[utf8]{inputenc}
\usepackage[english]{babel}
\usepackage{comment}
\numberwithin{equation}{section}

\usepackage{hyperref}

\usepackage[margin=1in]{geometry}
\usepackage{mathtools}
\usepackage{enumitem}
\usepackage{calc}

\usepackage{amsthm}
\usepackage{amssymb}
\usepackage{amsfonts}
\usepackage[mathcal]{eucal}
\usepackage{mathrsfs}
\usepackage{stix}

\swapnumbers
\theoremstyle{plain}
\newtheorem{theorem}{Theorem}

\newtheorem{proposition}[theorem]{Proposition}
\newtheorem{lemma}[theorem]{Lemma}

\theoremstyle{definition}

\newtheorem{definition}[theorem]{Definition}

\newtheorem{example}[theorem]{Example}

\newtheorem{remark}[theorem]{Remark}
\newtheorem{remarks}[theorem]{Remarks}

\theoremstyle{remark}

\newcommand{\abs}[1]{\left|{#1}\right|}
\newcommand{\nrm}[1]{\left\|{#1}\right\|}
\newcommand{\nrmB}{\nrm{\bullet}}
\newcommand{\NrmB}{\Nrm{\bullet}}
\newcommand{\nrmT}[1]{\nrm{#1}_{T}}

\newcommand{\nrmS}[1]{\nrm{#1}_{S}}
\newcommand{\nrmone}[1]{\nrm{#1}_{\ell_1}}
\newcommand{\nrminf}[1]{\nrm{#1}_{\ell_{\infty}}}

\newcommand{\fb}{\varphi_{\bullet}}

\newcommand{\tphi}{\tilde{\varphi}}

\newcommand{\cP}{\mathcal{P}}
\newcommand{\cS}{\mathcal{S}}
\newcommand{\cT}{\mathcal{T}}

\newcommand{\sS}{\mathscr{S}}
\newcommand{\sT}{\mathscr{T}}
\newcommand{\ThUmN}{\mathrm{Th}_\mathrm{UmN}}
\newcommand{\les}{\leqslant}

\newcommand{\fL}{\mathfrak{L}}
\newcommand{\LNseq}{L_{\mathrm{Nseq}}}
\newcommand{\LUmN}{L_{\mathrm{UmN}}}

\newcommand{\cM}{\mathcal{M}}
\newcommand{\tM}{\widetilde{\cM}}
\newcommand{\NN}{\mathbb{N}}

\newcommand{\bN}{\mathbf{N}}

\newcommand{\RR}{\mathbb{R}}

\newcommand{\sumsube}{\ensuremath{\sum^{\subseteq}}}
\newcommand{\sumnotsube}{\ensuremath{\sum^{\not\subseteq}}}

\newcommand{\ThTs}{\mathrm{Th}_{\cT}}
\newcommand{\ThSch}{\mathrm{Th}_{\sS}}
\newcommand{\cU}{\mathcal{U}}

\newcommand{\xb}{x_{\bullet}}

\newcommand{\Nrm}[1]{{\left\vert\kern-0.25ex\left\vert\kern-0.25ex\left\vert #1 
    \right\vert\kern-0.25ex\right\vert\kern-0.25ex\right\vert}}

\newcommand{\dotle}{\mathrel{\dot\le}}

\DeclareMathOperator{\supp}{supp}

\DeclareMathOperator{\NRM}{\boldsymbol{N\!orm}}

\begin{document}

\title
{The Non-Definability of the Spaces of Tsirelson and Schlumprecht}

\author[Casazza, Dueñez, Iovino]{Peter G. Casazza, Eduardo Dueñez \and José N. Iovino}

\address{Department of Mathematics\\
  University of Missouri\\
  Columbia, MO 65211\\
  U.S.A.}
\email{casazzap@missouri.edu}

\address{Department of Mathematics\\
  The University of Texas at San Antonio\\
  San Antonio, TX 78249\\
  U.S.A.}
\email{eduardo.duenez@utsa.edu}
\email{jose.iovino@utsa.edu}

\thanks{Casazza was supported by NSF DMS 1609760, NSF DMS 1725455, and ARO W911NF-16-1-0008}

\thanks{Dueñez and Iovino were supported by NSF DMS-1500615}



 \begin{abstract}
We prove the impossibility of finding explicit finitary definitions of the spaces of Tsirelson and Schlumprecht in continuous first-order logic.
\end{abstract}

\maketitle

\section*{Introduction}

For decades now, it has been an open question whether Tsirelson's space~\cite{Tsirelson:1974,Figiel-Johnson:1974} admits a finitary explicit definition~\cite{Gowers:1995},~\cite[Question~3, page 201]{Odell:2002}.
In this paper we use ideas of model theory to prove the non-existence of a finitary definition of Tsirelson's space. 
Here,``finitary'' means ``first-order''.
We use Keisler's framework of model theory for general (real-valued) structures~\cite{Keisler2020}, which provides the most general approach for the logical analysis of structures endowed with real-valued functions.

Our argument has two main ingredients:
a model-theoretic one and an analytic one. 
The analytic ingredient is a non-uniform convergence argument that relies on a construction involving fast-growing functions. 
The model-theoretic ingredient of the proof relies on Shelah's theory of stability~\cite{Shelah:1990}.
The aforementioned non-uniform convergence result allows us to deduce that the Tsirelson  does not have a finitary definition  (in the sense of  Gaifman-Shelah~\cite{Gaifman:1976,Shelah:1971a})  in the complete theory of structure $(c_{00},\|\cdot\|_n)_{n\in \mathbb{N}}$, where $\|\cdot\|_n$ is the $n$-the approximant of the Tsirelson norm. See Theorem~\ref{thm:nondef-Tsirelson}.


Our approach extends to other implicitly defined normed spaces, although the key analytic ingredient must be proved on a case-by-case basis.
To illustrate this, in the last section of the paper, we extend our non-definability result to Schlumprecht's space~\cite{Schlumprecht:1991}.

The only technical prerequisite of the paper is familiarity with model theory for  real-valued structures.
We use the recent framework of Keisler's general structures~\cite{Keisler2020};
however, our arguments can be easily translated to the well-known formalism of metric structures \emph{à la} Ben Yaacov-Usvyatsov~\cite{Ben-Yaacov-Usvyatsov:2010}
(see remarks~\ref{rem:master-norm}).

The paper is organized as follows:
In section~\ref{sec:Tsirelson}, we recall the definition of the Tsirelson space of Figiel and Johnson and prove that the convergence of the Tsirelson approximates towards the Tsirelson norm is non-uniform.
In section~\ref{sec:undef-Tsirelson}, we introduce the model-theoretic point of view and exhibit the link between definability and uniform convergence. 
Recalling the non-uniformity lemmas of section~\ref{sec:Tsirelson}, we deduce the non-definability of the Tsirelson norm. This is the main result of the paper.
In section~\ref{sec:Schlumprecht}, we refine the approach of section~\ref{sec:undef-Tsirelson} to deduce the non-definability of Schlumprecht's space. 

We are grateful to Bill Johnson for encouraging us to collaborate in this project.

This paper is dedicated to the memory of Boris Tsirelson.
We were honored when he wrote to us expressing his pleasure for the result presented here.

\section{The Tsirelson space}
\label{sec:Tsirelson}

\subsection{The Tsirelson space of Figiel and Johnson}
\label{sec:Tsirelson-defns}

Let us start this section by recalling the construction of the Tsirelson space of Figiel-Johnson~\cite{Figiel-Johnson:1974}.
(The reader is referred to~\cite{Casazza-Shura:1989} for a comprehensive treatment of Tsirelson-like spaces.)

If $E,F$ are finite non-empty subsets of $\NN =
\{1,2,3,\dots\}$, we write $E\le F$ if $\max E \le \min F$, and $E<F$ if $\max E < \min F$.
We will also write $n \le E$ for $\{n\}\le E$.

Let $c_{00}$ denote the space of all finite real scalar
sequences, and let $\{e_j\}_{j=1}^{\infty}$ be the canonical basis of
$c_{00}$.
For any $x = \sum_{n=1}^{\infty}a_ne_n$ in $c_{00}$ and any $1\le E
\subset \NN$, define 
\[ 
Ex = \sum_{n\in E}a_ne_n.
\]
 A collection of $m\ge 1$ nonempty sets $E_1,E_2,\dots,E_m$ of natural numbers is \emph{admissible} if $m\le E_1 < E_2 <\dots < E_m$. 
 Define a sequence of norms $\{\nrm{\cdot }_k\}_{k=0}^{\infty}$ on $c_{00}$ inductively as follows.
 For fixed $x \in c_{00}$, say $x=\sum_{n=1}^{\infty}a_ne_n$, let
\[ \nrm{x}_0 = \max_{n}|a_n|\]
and, for $k\ge 0$,
\[ 
\nrm{x}_{k+1}
= \max\left( \nrm{x}_k,\ 
  \max \left\{\frac{1}{2}\sum_{i=1}^m\nrm{E_ix}_k:
    1\le m\le E_1 < E_2 <\dots < E_m\right\}
  \right),
\]
where the innermost maximum is taken over the set of all admissible collections of any size $m\ge 1$.
 For $x = \sum_{n=1}^{\infty}a_ne_n$, we have
\[ 
\nrm{x}_k \le \nrm{x}_{k+1} \le \sum_{n=1}^{\infty}|a_n|=\nrmone{x};
\]
thus, each norm $\nrm{\cdot}_k$ is 1-Lipschitz with respect to $\nrmone{\cdot}$.

\begin{definition}
\label{D:Tsirelson}
The \emph{Tsirelson norm~$\nrmT{\cdot}$} on~$c_{00}$ is defined as
\[
  \nrmT{x}=\lim_{k\rightarrow \infty}\nrm{x}_k.\qquad
  \text{($x\in c_{00}$.)}
\]
\emph{Tsirelson's space}, denoted~$\cT$, is the norm-completion of $(c_{00}, \nrmT{\cdot})$.
(As customary, the extension of the Tsirelson norm to Tsirelson's space is still denoted $\nrmT{\cdot}$.)
The norm $\nrm{\cdot}_k$ will be called be the $k$-th \emph{iterate} in the construction of the Tsirelson norm.
\end{definition}

\subsection{Non-uniform convergence of the Tsirelson approximants.}
\label{sec:nonunif-Tsir}
As in section~\ref{sec:Tsirelson-defns} above, $\nrm{\cdot}_k$ will denote the $k$-th iterate in the construction of the Tsirelson norm on~$c_{00}$.

As a preliminary step, we need to construct certain functions of rapid growth.
Recall that $g^{(m)}$ denotes the $m$-fold iteration $g\circ\dots\circ g$ of a function~$g$ on a set~$X$, with the usual convention that $g^{(0)}$ is the identity map on~$X$.)

\begin{definition}\label{def:f-fast-growth}
  Let functions $f_1, f_2, \dots, f_k, \dots$ on natural numbers be as follows.
  For all $m\in\NN$, let $f_1(m) = m^2$, and recursively define $f_{k+1}(m) = f_k^{(m)}(m)$ for $k\in\NN$.
\end{definition}
Thus, $f_2(m) = m^{2^{m-1}}$, but $f_k$ for $k > 2$ admits no explicit algebraic expression.

Bellenot had found essentially the same rapid-growth functions in connection, not with the convergence of approximants norms, but with James' non-constructive proof that the standard basis $(e_n)$ of Tsirelson's space is not subsymmetric~\cite{James1964,Bellenot:1984}.%
\footnote{More specifically, Bellenot's theorem shows that $(e_{n_k})$ is not equivalent to the standard basis $(e_n)$ of Tsirelson's space if and only if $(n_k)$ asymptotically grows faster than $(f_k(k))$ with $f_k$ as in definition~\ref{def:f-fast-growth} below.
This answers a question of Casazza as to whether ``reasonable'' sequences $(n_k)$ yield $(e_{n_k})$ equivalent to~$(e_n)$.}

The support $\{n\in\NN : a_n\ne 0\}$ of an element $x = \sum_{n\in\NN}a_ne_n\in c_{00}$ will be denoted $\supp(x)$.
For $c < d\in\NN$, we write $[c,d)$ for $\{c,c+1,\dots,d-1\}$.

\begin{theorem}\label{thm:Tsirelson-bounds}
For $k\ge 1$, let $\nrm{\cdot}_k$ be the $k$-th iterate in the definition of the Tsirelson norm.
For every $n\ge 2$ there are vectors $\{x_i\}_{i=1}^n$ in $c_{00}$ such that:
\begin{enumerate}
\item $x_j$ is supported on a subset of $\left[f_k^{(j-1)}(n),2\sqrt{f_k^{(j)}(n)}\right)$,
  and the sum $x = \sum_{j=1}^n x_j$ is supported on a subset of $\left[n, 2\sqrt{f_{k+1}(n)}\right)$,
\item $\nrm{x_i}_k=\frac{1}{2}$ for $i=1,2,\ldots,n$,
\item $\nrm{x}_k\le 1$,
\item $\nrm{x}_{k+1}\ge \frac{n}{4}$, and
\item $\nrmone{x} \le 2^{k-1}n$.
\end{enumerate}
\end{theorem}

We will prove Theorem~\ref{thm:Tsirelson-bounds} by induction on~$k$.

\begin{proof}[Proof of Theorem~\ref{thm:Tsirelson-bounds} for $k=1$ (Base step)]
  Fix $n\ge 2$.
For $i=1,2,\ldots,n$, let $m_i = f_1^{(i-1)}(n) = n^{2^{i-1}}$, and
\[
  E_i = [m_i,2m_i) = \{m_i,m_i+1,\ldots,2m_i-1\}.
\]
 Define vectors $\{x_i\}_{i=1}^n$ by
\[
  x_i= \frac{1}{m_i}\sum_{j=m_i}^{2m_i-1}e_j.
\]
Clearly, $\supp(x_i) = E_i = [m_i, 2m_i) = \bigl[m_i, 2\sqrt{f_1(m)}\bigr)$ (since $f_1(m) = m^2\ge 2m$ for $m\ge 2$).
Thus, $x$ is supported on a subset of $\left[m_1, 2\sqrt{f_1(m_n)}\right) = \left[n, 2\sqrt{f_1^{(n)}(n)}\right) = \bigl[n, 2\sqrt{f_2(n)}\bigr)$, proving assertion~(1).
We also have $\nrmone{x} = \sum_{i=1}^n\nrmone{x_i} = \sum_{i=1}^n\sum_{j=m_i}^{2m_i-1}\frac{1}{m_i} = nm_i/m_i = n$, proving~(5).
For $i = 1,2,\ldots,n$, we have $\nrm{x_i}_0 = 1/m_i$, hence (2) follows from
\[
  \nrm{x_i}_1 = \frac{1}{2}\sum_{j=m_i}^{2m_i-1}\frac{1}{m_i}
  = \frac{m_i}{2m_i} = \frac{1}{2}.
\]
By disjointness of the supports of the~$x_i$, we have
\[
  \nrm{x}_2\ge \frac{1}{2}\sum_{i=1}^n
  \nrm{x_i}_1= \frac{1}{2}\sum_{i=1}^n \frac{1}{2}= \frac{n}{4},
\]
proving~(4).
It remains to prove~(3).
Since the norms $\nrm{x_i}_0 = 1/m_i$ decrease with~$i$, it follows that $\nrm{x}_0 = \nrm{x_1}_0 = 1/m_1 = 1/n$.
In order to find an upper bound for~$\nrm{x}_1$, we fix~$q$ and any admissible collection $q \le F_1< F_2 < \dots < F_q$ in order to find an upper bound for
\begin{equation*}
  \Sigma \coloneq  \frac{1}{2}\sum_{j=1}^q\nrm{F_jx}_0.
\end{equation*}
We may assume that all the sets $F_j$ are nonempty, and
\[ \bigcup_{j=1}^qF_j \subset \bigcup_{j=1}^nE_j,\]
since, otherwise, $F_j$ is unnecessarily capturing zero coefficients of~$x$.
Without loss of generality we may assume that $q\in F_1$.
(Otherwise, $\Sigma$ is possibly enlarged when one appends to $F_1$ the integers between $q$ and $\min F_1$.)
Henceforth, we fix $1\le p \le n$ such that $m_p \le q < m^2_p$.
There exists a unique index $1\le \ell \le q$ satisfying%
\footnote{We adopt the standard conventions $\bigcup_{j=b}^{b-1}Z_j  = \emptyset$ and $\sum_{j=b}^{b-1}z_j = 0$.}
\begin{itemize}
\item $\bigcup_{j=1}^{\ell-1}F_j \subset E_p$,
\item
$F_{\ell}\cap E_p \ne \emptyset$,\quad and
\item 
  $\bigcup_{j=\ell+1}^q F_j \subset \bigcup_{j=p+1}^nE_j$.
\end{itemize}
By the first and second properties above, $F_1,F_2,\dots,F_{\ell-1}$ and $F_{\ell}\cap E_i$ are $\ell$ nonempty disjoint subsets of~$E_p$, so $m_p = \#E_p \ge \ell$.
The three properties above imply corresponding inequalities:
\begin{itemize}
\item
$\nrm{F_jx}_0= \frac{1}{m_p}$, for $j = 1,2,\ldots,\ell$,\quad and
\item
$\nrm{F_jx}_0\le \frac{1}{m_{p+1}}$, for $j=\ell+1,\ldots,q$.
\end{itemize}
For convenience, let $m_{n+1} = +\infty$.
We have
\begin{align*}
  \Sigma
  &= \frac{1}{2}\left [
    \sum_{j=1}^{\ell-1}\nrm{F_jx}_0
    + \nrm{F_{\ell}x}_0
    + \sum_{j=\ell+1}^q\nrm{F_jx}_0
    \right ]\\ 
  &\le \frac{1}{2}
    \left [
    \frac{\ell-1}{m_p}
    + \frac{1}{m_p}
    + \sum_{j=\ell+1}^q
    \frac{1}{m_{p+1}}\right ]\\ 
  &= \frac{1}{2}
    \left[
    \frac{\ell}{m_p}+\frac{q-\ell}{m_{p+1}}
    \right]\\
  &\le \frac{1}{2}
    \left[
    \frac{\ell}{m_p}+\frac{m_p^2}{m_{p+1}}
    \right]
    \qquad\text{(since $\ell>0$ and $q < m_p^2$)}\\
  &\text{(the term $m_p^2/m_{p+1}$ above is zero when $\ell=q$)}\\
  &\le \frac{1}{2}\left[1 + 1\right]
  \qquad\text{(since $\ell\le m_p$ and $m_{i+1} \ge m_i^2$)}\\
  &= 1. 
\end{align*}
This proves~(3), finishing the proof of Theorem~\ref{thm:Tsirelson-bounds} for $k=1$.
\end{proof}

\begin{proof}[Proof of Theorem~\ref{thm:Tsirelson-bounds} (Inductive step)]
Assume the statement in the Theorem holds for some fixed $k\ge 1$.
Fix $n\ge 2$.
For $1\le i\le n$, let $m_i = f_{k+1}^{(i-1)}(n)$, and for $1\le j\le m_i$, let $m_{ij} = f_k^{(j-1)}(m_i)$.
Apply the inductive assumption with $m_i$ in place of~$n$ to obtain $m_i$~vectors, say $\{x_{ij}\}_{j=1}^{m_i}$, such that:
\begin{enumerate}
\item $\supp(x_{ij}) \subset \bigl[m_{ij}, 2\sqrt{f_k(m_{ij})}\bigr)$ for $1\le j\le m_i$, and the sum $x_i = \sum_{j=1}^n x_{ij}$ is supported on a subset of $\left[m_i, 2\sqrt{f_{k+1}(m_i)}\right)$;
\item $\nrm{x_{ij}}_k=\frac{1}{2}$ for $1\le j\le m_i$;
\item $\nrm{x_i}_k\le 1$,
\item $\nrm{x_i}_{k+1}\ge \frac{m_i}{4}$, and
\item $\nrmone{x_i} \le 2^{k-1}m_i$.
\end{enumerate}

For each~$i$, let
\[
  y_i = \frac{x_i}{2\nrm{x_i}_{k+1}}.
\]
As shown above, the elements $y_i$ ($1\le i\le n$) together with their sum $y = \sum_{i=1}^n y_i$ satisfy~(1) of Theorem~\ref{thm:Tsirelson-bounds} with $k+1$ in place of~$k$ ($y$ is supported on a subset of $\bigl[n, 2\sqrt{f_{k+1}(m_n)}\bigr)$, and $f_{k+1}(m_n) = f_{k+1}^{(n)}(n) = f_{k+2}(n)$).
 Next, we have
 \begin{align*}
   \nrm{y_i}_{k+1} &= \frac{1}{2},\\
   \nrm{y_i}_k &= \frac{\nrm{x_i}_k}{2\nrm{x_i}_{k+1}}
    \le \frac{1}{2m_i/4} = \frac{2}{m_i},\\
   \nrm{y}_{k+2} &\ge \frac{1}{2}\sum_{i=1}^n\nrm{y_i}_{k+1} =\frac{1}{2}\sum_{i=1}^n\frac{1}{2} = \frac{n}{4},
 \end{align*}
 proving properties~(2) and~(4) for~$k+1$.
Property~(5) follows from
\begin{align*}
  \nrmone{y} = \sum_{i=1}^n \nrmone{y_i}
  = \frac{1}{2} \sum_{i=1}^n \frac{\nrmone{x_i}}{\nrm{x_i}_{k+1}}
  \le \frac{1}{2} \sum_{i=1}^n \frac{2^{k-1}m_i}{m_i/4}
  = 2^kn.
\end{align*}

It remains to show that property~$(3)$ holds, i.e., that $\nrm{y}_{k+1}\le 1$.
Fix~$q$ and any admissible collection $q \le F_1< F_2 < \dots < F_q$.
As in the proof of the case~$k=1$, we may assume that $q = \min F_1$, and that for some~$p$ we have $m_p\le q < 2\sqrt{f_{k+1}(m_p)}$.
For each $i,j$, let $E_{ij}$ be the support of~$y_{ij}$, and let $E_i = \bigcup_{j=1}^{m_i}E_{ij}$ be the support of~$y_i$. 
There exists a unique~$\ell$ ($1\le\ell\le n$) such that (\emph{i})~$\bigcup_{i=1}^{\ell-1}F_i\subset E_p$, (\emph{ii})~$F_{\ell}\cap E_p\ne \emptyset$, and (\emph{iii})~if $\ell<n$, then $F_{\ell+1} \cap E_p=\emptyset$.
For $1\le i\le n$ and $1\le j\le q$, let $z_{ij} = F_jy_i = (E_i\cap F_j)y$.
We have $z_{ij}=0$ if $i<p$ (by~(\emph{i}) above), and $z_{pj}=0$ if $j>\ell$ (by~(\emph{iii})).
We seek an upper bound for
\begin{equation*}
  \Sigma \coloneq  \frac{1}{2}\sum_{j=1}^q\nrm{F_jy}_k
  = \frac{1}{2}\sum_{j=1}^q\nrm{\sum_{i=p}^nz_{ij}}_k.
\end{equation*}
By the triangle inequality (and omitting zero terms):
\begin{align*}
  \Sigma \le \frac{1}{2}\sum_{i=p}^n\sum_{j=1}^q\nrm{z_{ij}}_k
    &= \frac{1}{2} \sum_{j=1}^{\ell}\nrm{z_{pj}}_k +
      \frac{1}{2} \sum_{i=p+1}^n\sum_{j=\ell}^q\nrm{z_{ij}}_k \\
  &= \Sigma' + \Sigma''.
\end{align*}
On the one hand,
\[
  \Sigma' \coloneq  \frac{1}{2}\sum_{j=1}^{\ell}\nrm{z_{pj}}_k\le \nrm{y_p}_{k+1}
  = \frac{1}{2}.
\]
On the other hand, assuming $p < n$ (otherwise the quantity~$\Sigma''$ is zero), since $z_{ij} = F_jy_i$ by definition and $v\mapsto F_jv$ is norm-decreasing, we have
\begin{equation*}
 \Sigma''
  \coloneq  \frac{1}{2}\sum_{i=p+1}^n\sum_{j=\ell}^q\nrm{z_{ij}}_k
\le \frac{1}{2}\sum_{i=p+1}^n\sum_{j=\ell}^q\nrm{y_i}_k.
\end{equation*}
Under the temporary assumption $n\ge 3$, since $k\ge 1$, we have
\begin{equation*}
f_{k+1}(n) \ge f_2(n) = f_1^{(n)}(n) \ge f_1^{(3)}(n) = n^8.
\end{equation*}
It follows that
\begin{align*}
  \Sigma''
  &\le \frac{1}{2}nq\max_{p+1\le i\le n}\nrm{y_i}_k
    \le \frac{nq}{m_{p+1}}
\qquad(\nrm{y_i}_k\le 2/m_i \le 2/m_{p+1})\\
&\le \frac{2n\sqrt{f_{k+1}(m_p)}}{m_{p+1}}
\qquad\left(q < 2\sqrt{f_{k+1}(m_p)}\right)\\
  &\le \frac{2n}{\sqrt{f_{k+1}(n)}} \le \frac{2n}{\sqrt{n^8}}
    = \frac{2}{n^3} \le \frac{2}{3^3}
    < \frac{1}{2}
\qquad(m_{p+1} \ge f_{k+1}(m_p) \ge f_{k+1}(n) \ge n^8).
\end{align*}
When $n=2$, since $k\ge 1$, we have $m_1=2$ and $m_2 = f_{k+1}(m_1) \ge f_2(2) = 16$.
If $\Sigma''\ne 0$, then the range of index~$i$ in the outer sum defining~$\Sigma''$ is simply $i=2=n=p+1$, so we may replace $nq$ by $q$ on the first line of the estimate above:
\begin{equation*}
  \Sigma'' \le \frac{q}{m_2} < \frac{2\sqrt{m_2}}{m_2}
  = \frac{2}{\sqrt{m_2}} \le \frac{2}{\sqrt{16}} = \frac{1}{2},
\end{equation*}
since $m_2 = f_{k+1}(m_1) = f_{k+1}(2) = f_k^{(2)}(2) \ge f_1^{(2)}(2) = 16$.

It follows that
\[ \frac{1}{2}\sum_{i=1}^q\nrm{F_iy}_k \le 1,\]
and hence $\nrm{y}_{k+1}\le 1$.
This proves property~(3) for $k+1$, completing the inductive step, and the proof of Theorem~\ref{thm:Tsirelson-bounds}.
\end{proof}

\begin{proposition}\label{prop:triangle-Tsirelson}
  There exist sequences
  \begin{itemize}
  \item $(x_i)_{i=1}^{\infty}$ in $c_{00}$, and 
  \item $(k_j)_{j=1}^{\infty}$ strictly increasing in~$\NN$, 
  \end{itemize}
  such that
  \begin{itemize}
  \item $\nrm{x_i}_{k_j} \le 1/3$\quad for $j\le i$,
  \item $\nrm{x_i}_{k_j} \ge 2/3$\quad for $j > i$, and
  \item $\lim_{j\to\infty}\nrm{x_i}_{k_j} = \nrmT{x_i} = 1$\quad for all~$i$.
  \end{itemize}
\end{proposition}
The proof of Proposition~\ref{prop:triangle-Tsirelson} hinges on the following Lemma.

\begin{lemma}\label{lem:Tsirelson-gaps}
  Given $k\in\NN$, there exist $k'>k$ and $x\in c_{00}$ such that
  \begin{align*}
    \nrm{x}_k &\le 1/3,&
                         \nrm{x}_{k'} &\ge 2/3,&
                                        \quad\nrmT{x} &= 1.
  \end{align*}  
\end{lemma}
\begin{proof}
  By Theorem~\ref{thm:Tsirelson-bounds} with $n=12$, given $k$ there exists $y\in c_{00}$ such that $\nrm{y}_k\le 1$ and $\nrm{y}_{k+1}\ge 3$.
  Since $N \coloneq  \nrmT{y} = \sup_l\nrm{y}_l \ge \nrm{y}_{k+1} = 3$, there exists $k'$ with $\nrm{y}_{k'}\ge \max\bigl\{\frac{2}{3}N, 3\bigr\}$.
  The vector $x = y/N$ satisfies:
  \begin{itemize}
  \item $\nrm{x}_k  = \nrm{y}_k/N \le 1/N \le 1/3$,
  \item $\nrm{x}_{k'} = \nrm{y}_{k'}/N \ge \frac{2}{3}N/N = 2/3$,\quad and
  \item $\nrmT{x} = \nrmT{y}/N = N/N = 1$.
  \end{itemize}
  Since the sequence $(\nrm{x}_k)_{k=1}^{\infty}$ is nondecreasing, the first two properties above imply that $k'>k$.
\end{proof}

\begin{proof}[Proof of Proposition~\ref{prop:triangle-Tsirelson}]
  Given $k\in\NN$, let $x_k \coloneq  x$ and $f(k) \coloneq  k'$ be those given by Lemma~\ref{lem:Tsirelson-gaps}.
  Define a sequence $(k_j)_{j=1}^{\infty}$ in $\NN$ recursively by
  \begin{itemize}
  \item $k_1 = 1$, \quad and
  \item $k_{j+1} = f(k_j)$ \quad for $j\in\NN$,
  \end{itemize}
  i.e., $k_j = f^{(j-1)}(1)$.
  We claim that the sequences $(x_{k_i})_{i=1}^{\infty}$ and $(k_j)_{j=1}^{\infty}$ satisfy the properties stated in Proposition~\ref{prop:triangle-Tsirelson}.
  Since $f$ is strictly increasing per Lemma~\ref{lem:Tsirelson-gaps}, so is the sequence~$(k_j)$;
  in particular, $k_j\to\infty$ as $j\to\infty$.
  Write $N(i,j)$ for $\nrm{x_{k_i}}_{k_j}$. 
  By Lemma~\ref{lem:Tsirelson-gaps}, the monotonicity of the family $\{\nrm{\cdot}_k\}$ and of the sequence $(k_j)$, we have
  \begin{itemize}
  \item For $j\le i$: $N(i,j) \le N(i,i) \le 1/3$,
  \item For $j > i$: $N(i,j) \ge N(i,i+1) \ge 2/3$,
  \item For all $x$: $\lim_{j\to\infty}\nrm{x}_{k_j} = \nrmT{x}$\quad(since $k_j\to\infty$ as $j\to\infty$),
  and
  \item For all $i$: $\nrmT{x_i} = 1$.\qedhere
  \end{itemize}
\end{proof}

\begin{proposition}\label{prop:double-limit-Tsirelson}
  There exist sequences $(x_i)_{i=1}^{\infty}$ in $c_{00}$, and $(k_j)_{j=1}^{\infty}$ strictly increasing in~$\NN$, such that
  \begin{itemize}
  \item the limit $\lambda_j \coloneq  \lim_{i\to\infty}\nrm{x_i}_{k_j}$ exists for each $j\in\NN$, and
  \item the limit $\Lambda \coloneq  \lim_{j\to\infty}\lambda_j$ exists, satisfies $\Lambda \le 1/3$, and
  \item $\lim_{j\to\infty}\nrm{x_i}_{k_j} = \nrmT{x_i} = 1$
    for all~$i\in\NN$.
  \end{itemize}
\end{proposition}
\begin{proof}
  Let $(y_i)_{i=1}^{\infty}$ in~$c_{00}$ and $(l_j)_{j=1}^{\infty}$ be sequences satisfying the conclusions of Proposition~\ref{prop:triangle-Tsirelson}.
  Let $N(i,j) \coloneq  \nrm{y_i}_{l_j}$.
  For fixed~$j$, the sequence~$\{N(\cdot,j)\}$ eventually takes values in~$[0,1/3]$;
  by sequential compactness of~$[0,1/3]$, it has a convergent subsequence, say $\{N(n_i,j)\}_{i=1}^{\infty}$.
  The choice of indexes $(n_i)_{i=1}^{\infty}$ realizing the subsequence depends on the fixed choice of~$j$;
  we will write $n_{ij}$ instead of $n_i$ in order to exhibit the dependence explicitly.
  Thus, $\lambda_j \coloneq  \lim_{i\to\infty}N(n_{ij},j)$ exists for each~$j$, and $\lambda_j\le 1/3$.
  Let $x_i = y_{n_{ii}}$ for each~$i$, so $(x_i)$ is the ``diagonal'' subsequence of the family of sequences $\{(y_{n_{ij}})_i:j\in\NN\}$.
  Once more, by sequential compactness, there is a subsequence $(l_{m_j})_{j=1}^{\infty}$ of $(l_j)$ such that the limit $\Lambda \coloneq  \lim_{j\to\infty}\lambda_j$ exists, and necessarily $\Lambda\le 1/3$.
  The sequences $(x_i)$ and $(k_j)$ with $k_j = l_{m_j}$ satisfy the following properties:
  For fixed $j\in\NN$, $(n_{ii})_i$ is a subsequence of~$(n_{ij})_i$, so
  \begin{equation*}
    \lambda_j
    = \lim_{i\to\infty}N(n_{ij},j)
    = \lim_{i\to\infty}N(n_{ii},j)
    = \lim_{i\to\infty}\nrm{x_i}_{k_j}, 
 \end{equation*}
  where the latter (subsequential) limit necessarily exists because $(x_i)_{i=1}^{\infty}$ is a subsequence of~$(y_{n_{ij}})_{i=1}^{\infty}$.
  Finally, since $(k_j)$ is a subsequence of~$(l_j)$,
 \begin{equation*}
   \lim_{j\to\infty}\nrm{x_i}_{k_j}
   = \lim_{j\to\infty}\nrm{y_{n_{ii}}}_{l_j}
   = \nrm{y_{n_{ii}}}_T
   = 1.
   \qedhere
 \end{equation*}
\end{proof}

\section{First-Order Non-Definability of the Tsirelson Space}
\label{sec:undef-Tsirelson}

\subsection{Uniformly Multinormed Structures}
\label{sec:UmN}
Let $V$ be any real vector space endowed with a pointwise bounded collection $\bN$ of seminorms, i.e., such that for each $x\in V$ there exists $C\ge 0$ such that $N(x)\le C$ for all $N\in\bN$.
The \emph{(concrete) uniformly multi-normed (UmN) structure associated to $(V,\bN)$} is the structure
\[
cM = \langle \bN, V_n, 0_n, +_n, \cdot_{t,n}, \iota_{m,n}, \NRM_k: m,n\in\NN,\ m\le n,\ t\in\RR \rangle,
\]
where $\bN$ and $V_n$ (for $\in\NN$) are distinct sorts of $\cM$ and  for all $m\le n\in\NN$:
\begin{itemize}
\item $V_n = \{x\in V: \text{$N(x)\le n$ for all $N\in\bN$}\}$;
\item $0_n$ is the zero of~$V$,
\item $\iota_{m,n}$ is the inclusion $V_m\hookrightarrow V_n$;
\item $+_n$ and $\cdot_{t,n}$ are (the restrictions of) the sum and scalar product-by-$t$ on~$V$ to functions $V_n\times V_n\to V_{2n}$ and $V_n\to V_l$ where $l = \lceil \left|tn\right|\rceil$ (the least integer $\ge\left|tn\right|$).
\item $\NRM_n(N,x) = N(x)$ for each $N\in \bN$ and $x\in V_n$.
\end{itemize}
  The \emph{language for UmN structures} is the first-order (real-valued) language $\LUmN$ for general structures whose vocabulary consists of the sorts, constants, functions, and the predicates $\NRM_n$ above.

We shall allow a slightly more general notion of concrete UmN structure, so as to obtain a first-order axiomatizable class (up to canonical identifications%
\footnote{See Remarks~\ref{rem:canonical-UmN-identifications}}) of general structures in Keisler's sense.
We shall only require that $V_1 = \bigcap_{N\in\bN}B_N(r_N)$ where, for each $N\in\bN$, we have $0 < r_N\le 1$ and $B_N(r_N)$ is either the open or closed ball of radius $r_N$ relative to the seminorm~$N$ (some balls may be open and others closed)%
\footnote{The sorts $V_n$ first described above are obtained when the balls $B_N(1)$ are all of radius $r_N=1$ and closed.}
and, accordingly, letting $V_n \coloneq nV_1$ for all $n\in\NN$.
However, in this more general setting, $V'\coloneq\bigcup_nV_n$ may be a proper subset of~$V$ (in which case, nevertheless, one obtains exactly the same UmN structure from $(V', \bN\restriction V')$).
Using this more general notion of concrete UmN structure, a given nontrivial $(V,\bN)$ typically admits many formally different associated UmN structures (i.e., in different Keisler isomorphism classes).

The first-order theory (in real-valued logic) $\LUmN$-theory of concrete UmN structures is denoted~$\ThUmN$.
Its models are \emph{(abstract) uniformly multinormed (UmN) structures}, and consist of the following:
\begin{itemize}
\item A set (sort)~$\bN$, called the \emph{norms sort}.
    \item For all natural numbers $m,n$:
      \begin{itemize}
      \item a set (\emph{``sort''}) $V_{n}$ and an element $0_{n}\in V_n$;
      \item a real-valued predicate $\NRM_{n}$ on $\bN\times V_{n}$;
      \item for every $t\in\RR$, a function $\cdot_{t,n}: V_{n}\to V_{\lceil tn\rceil}$;
      \item an operation $+_{n,m} : V_{n}\times V_{m} \to V_{n+m}$;
      \item if $n\le m$, a function $\iota_{n,m}: V_{n}\to V_{m}$;
      \end{itemize}
    \end{itemize}
    and are models of~$\ThUmN$.
    Indeed, the zero elements and operations implicitly define a vector space~$V$ and a collection $\bN$ of seminorms therein (one may construct these explicitly from the direct limit $V\coloneq\varinjlim V_n$---relative to the identifications $\iota_{m,n}$---with operations $+$ and $\cdot$ induced by the operations $+_m$ and $\cdot_{t,n}$, plus one seminorm corresponding to each element $N\in\bN$---it is induced by the functions $\NRM_n(N,\cdot)$ as $n$ varies).
The theory of UmN structures is necessary and sufficient to ensure that $V_1$ is the intersection of open or closed balls (in the seminorms $N\in\bN$) of radii at most~$1$.


\begin{remarks}\label{rem:canonical-UmN-identifications}
  \begin{enumerate}
  \item Whenever sensible, if $\cM$ is an UmN structure, $N$ an element of $\bN^{\cM}$, and $x\in V_n^{\cM}$, we write $N(x)$ as an alias for the syntactically clumsy $\NRM_n^{\cM}(N, x)$.
  \item If $\cM$ is a concrete UmN structure obtained, so is its \emph{reduction} $\cM_{\mathrm{Red}}\coloneq\cM/\!\equiv$ obtained under identification of elements of each $V_n^{\cM}$ by the equivalence relation
    \begin{equation*}
      x\equiv y\quad\Leftrightarrow\quad\text{$N(y-x)=0$ for all $N\in\bN^{\cM}$},
    \end{equation*}
    and the correspondingly induced operations and interpretations of the Norm predicates.
    Moreover, $\cM$ and $\cM_{\mathrm{Red}}$ isomorphic in Keisler's sense, although we prefer to say that they are in the same isomorphism class (one of the peculiarities of the theory of general classes is that there need not exist a bijective isomorphism between the universes of $\cM$ of $\cM_{\mathrm{Red}}$ ---the quotient maps $a\mapsto a/\!\equiv$ from the universe of $\cM$ onto the universe of $\cM_{\mathrm{Red}}$ need not be injective).
  \end{enumerate}
\end{remarks}
In summary, every UmN structure is concrete, up to reduction and isomorphism.

Although the technical definition of UmN structures is needed to apply the Keisler formalism to spaces endowed with multiple norms, the preceding remarks essentially allow abstracting the sorts $V_n$ and regard a uniformly multinormed structure~$\cM$ as a ``metastructure'' $\cM^{*} = \langle V, 0, +, \cdot, \bN, \NRM\rangle$ with a single vector sort $V = \bigcup_nV_n$ with zero~$0$, operations $+$ and~$\cdot$, and $\bN$ a set of formal names (labels) for seminorms on~$\bN$ realized via the single predicate $\NRM$ (the latter ingredients of $\cM^{*}$ are naturally induced from those of~$\cM$).
A critical feature of the metapredicate $\NRM^{\cM}$ is its local boundedness on~$V^{\cM}$ (actually, $\NRM^{\cM}(N,\cdot)\le 1$ on~$V_1^{\cM}$ for all $N\in\bN^{\cM}$) as alluded by the adverb ``uniformly'' (although $\NRM^{\cM}(N,\cdot)$ is typically unbounded on all of~$V^{\cM}$).

In what follows, we mostly pretend away the technicalities of the multisorts $V_n$, and treat UmN metastructures in the above sense as though they are \emph{bona fide} Keisler general structures.
In particular, we shall treat UmN structures as concrete ones (with $V_k\subseteq V_l$ for $k\le l$, in particular).
The phenomenon that many different (say, concrete) UmN structures may yield the same metastructure hints at an external notion of isomorphism classes (i.e., isomorphism classes of $\equiv$-reduced UmN \emph{meta}structures)
typically larger than Keisler isomorphism classes.
However, we shall not pursue this avenue of research presently as it is not relevant to the considerations of this manuscript.
Still, when the multisorted nature of UmN structures is obscured by the metastructural viewpoint, we shall comment accordingly.

\begin{remarks}\label{rem:seminorms-internal-extension}
  \begin{itemize}
  \item The peculiar sort~$\bN$ whose elements formally name seminorms on~$V$ (rather than the seemingly more natural viewpoint of regarding each individual seminorm as a predicate) allows us to use the theory of stability and definability of types to study norms regarded as objects in the universe of the (meta)structure.
  \item Given a UmN structure~$\cM$ (say, concrete, at least for purposes of exposition) and any collection $S$ of seminorms on~$V^{\cM}$ such that $N(x)\le 1$ for all $N\in S$ and $x\in V_1^{\cM}$, one obtains an elementary extension $\tM\succeq\cM$ by keeping the same vector sort $V^{\tM}\coloneq V^{\cM}$, letting $\bN^{\tM} \coloneq S\sqcup\bN^{\cM}$, and extending $\NRM^{\cM}$ to $\bN^{\tM}\times V^{\tM}$ by $\NRM^{\cM}(N,x)\coloneq N(x)$ for $N\in S$ and $x\in V^{\tM}$.
  \end{itemize}
\end{remarks}

For the rest of this section, ${\cM}$ will be a fixed multi-normed structure associated to $(V,\bN)$ (thus $V$ is a real vector space and and $\bN$ is a pointwise bounded collection of seminorms on $V$) and we will maintain all the notational conventions introduced in this subsection.

\subsection{Types and Definability in UmN Structures}
\label{sec:types}


For our applications, we shall focus exclusively on \emph{mixed-sort types} in UmN structures; these types describe potential properties of object(s) of one sort ($V$ vs.\ $\bN$) in terms of parameters in the other sort.

  Fix a formula $\varphi(\overline{N};\overline{z})$ on an $m$-tuple $\overline{N}$ of formal variables of sort~$\bN$, and an $n$-tuple~$\overline{z}$ of variables of sort~$V$.
  Given an UmN-structure~$\cM$ and a set~$A\subseteq V^{\cM}$, the \emph{$\varphi$-type of~$\overline{N}\in(\bN^{\cM})^m$ (with parameters in~$A$)} is the function
  \begin{equation}\label{eq:norm-type}
    \begin{split}
      \varphi_N^{\cM}:
      A^n &\to \RR \\
      (x_1,\dots,x_n)
      &\mapsto
      \varphi^{\cM}(\overline{N}; x_1, \dots, x_n).
    \end{split}
  \end{equation}
  

  The notion of~$\tphi$-types, which are dual to the $\varphi$-types above, is obtained by formally exchanging the roles of the variables $\overline{N}$ and $\overline{x}$, i.e., $\tphi$ types are types for the formula $\tphi(\overline{x};\overline{N})\equiv\varphi(\overline{N};\overline{x})$.
Explicitly, the \emph{$\tphi$-type} of the $n$-tuple~$\overline{x}$ of elements of the universe~$V^{\cM}$ of an $\LNseq$-structure~$\cM$, with parameters in $B\subseteq\bN^{\cM}$, is the function
\begin{equation*}
  \begin{split}
    \tphi_{\overline{x}} :
    B^m&\to \RR\\
    \overline{N} &\mapsto \varphi^{\cM}(\overline{N}; \overline{x}).
  \end{split}
\end{equation*}
The $\varphi$-type of a tuple~$\overline{N}\in (\bN^{\cM})^m$ is said to be \emph{realized} in~$\cM$.
The collection of these is denoted $\cP_{\varphi}^{\cM}(A)$, where $A\subseteq V^{\cM}$ is the set of parameters of these types.
The collection $\cS_{\varphi}(A)$ of all $\varphi$-types (with parameters in~$A$) consists of the types of all $m$-tuples $\overline{N}$ of elements of sort $\bN$ in all elementary extensions $\langle\tM, a: a\in A\rangle \succeq \langle\cM, a:a\in A\rangle$.
The set $\cS_{\varphi}(A)$ inherits the subspace topology from the product space~$R^{A^m}$.
The formally multisorted nature of (the metasort) $V$ implies that $A\subseteq V_k^{\cM}$ for some $k\in\NN$, and furthermore $\cS_{\varphi}(A)$ is compact (by logical compactness).%
\footnote{Formally, the variables in the tuple $\overline{x}$ are necessarily of sort $V_k$ for some~$k$, so $\varphi$ \emph{sensu stricti} only involves parameters on $A\preceq V_k^{\cM}$ for (any) such~$k$.
On the other hand, since sort $\bN$ is discrete by definition, the parameter set $B\subseteq\bN^{\cM}$ is otherwise completely arbitrary.}
These considerations translate \emph{mutatis mutandis} to the sets $\cS_{\tphi}(B)$ of all dual types (realized in a fixed model~$\cM$, in the case of~$\cP^{\cM}_{\tphi}(B)$).%
\footnote{The spaces of dual types \emph{a priori} each consist of types $\tphi_{\overline{x}}$ with $\overline{x}$ of sort $V^k$ for some~$k$ depending only on~$\tphi$ itself.}


\subsubsection{Semidefinable Global Predicates and Definable Types.}
\label{sec:definable-types}

Fix a formula~$\varphi(\overline{N};\overline{x})$.
Throughout this section~$\cM$ is an arbitrary $\LUmN$-structure and $(\overline{N}_k)$ is a fixed sequence of $m$-tuples in~$\bN^{\cM}$.
The \emph{global $\tphi$-predicate~$P$ of kind~$V^n$ semidefined by the ultrafilter~$\cU\in\beta\NN$ and~$(\overline{N}_k)$} is the class of all real-valued functions
\begin{equation*}
  \begin{split}
    P^{\tM}: (V^{\tM})^n &\to\RR\\
    \overline{x} &\mapsto \lim_{\cU,k}\varphi^{\tM}(\overline{N}_k;\overline{x})
  \end{split}
\end{equation*}
as $\tM$ ranges over elementary extension of~$\cM$.          
Concretely, each $P^{\tM}$ is obtained in all $\tM\succeq\cM$ as the \emph{same} ultralimit of the sequence $(\varphi_{N_k}^{\tM})$ of types realized in $\cM$.


  A semi-definable global $\tphi$-predicate~$P$ is \emph{definable over a set of parameters~$B\subseteq\bN^{\cM}$} if $P^{\tM}$ is equal to the uniform limit of some sequence of types (with parameters in $V^{\tM}$) that are realized in~$\cM$ (though not necessarily of $\varphi$-types).
More concretely, $P$ is indeed so definable precisely if each $P^{\tM}$ is obtained as the same continuous combination of the collection $\langle\varphi_{\overline{N}}^{\tM} : \overline{N}\in B^m\rangle$, i.e., if there exists a continuous function~$C : \RR^{B^m}\to\RR$ such that
\begin{equation*}
  P^{\tM}(\overline{x})
  = C\bigl(\langle\varphi(\overline{N};\overline{x})^{\cM}:
  \overline{N}\in B^m\rangle\bigr)
  \quad
  \text{for all $\overline{x}\in (V^{\tM})^n$.}
\end{equation*}

A type $p\in \cS_{\varphi}^{\cM}$ is \emph{(semi)definable over a set of parameters~$B$} if there exists a (semi)definable global $\tphi$-predicate $P_p$ over~$B$ such that
\begin{equation}
  \label{eq:definition-scheme}
  p(\overline{x}) = P^{\cM}_p(\overline{x})
  \quad\text{for all $\overline{x}\in (V^{\cM})^n$.}
\end{equation}
Such predicate $P_p$ \emph{defines~$p$}.

\begin{remark}\label{rem:semidef-predic-type}
  %
  A definable $\tphi$-predicate $P$ over~$\cM$ admits an extension to a continuous function~$\cS_{\tphi}^{\cM}\to\RR$, which we also call~$P$ by an abuse of notation.
  Via its definition~$P_p$, a definable $\varphi$-type $p$ over~$\cM$ admits a natural extension to a continuous function on~$\cS_{\tphi}^{\cM}$.
  For details on definability and stability in the context of real-valued logic, the reader is referred to the literature~\cite{Ben-Yaacov-Usvyatsov:2010,BenYaacov2014,Keisler2020}.
\end{remark}

\subsubsection{Norm-Sequence Structures}
\label{sec:Nseq}

The \emph{language $\LNseq$ for norm sequences} expands $\LUmN$ with constant symbols (of the \emph{approximant norms}) $\nrmB_k$ ($k\in\NN$) and $\NrmB$ (the symbol for the \emph{master norm}) of sort~$\bN$.

We shall use the syntactic aliases $\nrm{x}_k$ (resp., $\Nrm{x}$) for $\NRM(\nrmB_k,x)$ (resp., for $\NRM(\NrmB,x)$).
The function $x\mapsto\NRM^{\cM}(\nrmB_k^{\cM},x)$ on $V^{\cM}$ will be denoted $\nrm{\cdot}^{\cM}_k$ (or just $\nrm{\cdot}_k$, if $\cM$ is clear), and similarly for $\Nrm{\cdot}$ ($=\Nrm{\cdot}^{\cM}$).

  
Given an Nseq structure~$\cM$ and an $n$-tuple $\overline{x}$ elements of the vector sort of~$\cM$, by formally identifying $\nrmB_k$ with its index~$k$, the $\tphi$-type $\tphi_{\overline{x}}$ ---with parameters in the collection $\bN_{\NN}^{\cM}\coloneq\{\nrmB_k^{\cM}: k\in\NN\}$ of approximant norms---is a bounded real sequence~$(t_k)_{k\in\NN}$, $t_k = \varphi^{\cM}(\nrmB_k;\overline{x})$.
Accordingly, we may regard $\cS^{\cM}_{\tphi}(\bN^{\cM}_{\NN})$ as a topological subspace of~$\RR^{\NN}$ (with the product topology, i.e., the topology of pointwise convergence).

\begin{remarks}\label{rem:master-norm}
  \emph An Nseq structure may have additional functions, constants, predicate interpretations, and possibly even other sorts.
  However, as long as the vocabulary~$L\supseteq\LNseq$ for such a structure is fixed \emph{and countable,} Keisler's Expansion Theorem for general such $L$-structures implies that any $L$-theory $T$ has a pre-metric expansion with a pseudo-metric approximate distance~$\mathrm{d}$.
  Although $\mathrm{d}$ is not canonical, the metric topology induced in models of~$T$ is.
  Thus, models of such a theory~$T$ are (multisorted) metric spaces, up to canonical reduction and metric equivalence.
  
  For structures with the exact vocabulary~$\LNseq$ and any $\LNseq$-theory $T$ including the sentences $\sup_{x\in V_k}\NRM_k(\nrmB_n,x) \dotle \sup_{x\in V_k}\NRM_k(\NrmB,x)$ for $k,n\in\NN$, one may regard the sort $\bN$ as a discrete space, and the sorts $V_k$ as metrized by $\mathrm{d}(x,y)\coloneq\Nrm{y-x}$ in models of~$T$.
  Thus, for Tsirelson and Schlumprecht structures as defined below, the formalism of metric structures \emph{à la} Ben Yaacov-Berenstein-Henson-Usvyatsov~\cite{Ben-Yaacov-Usvyatsov:2010,Ben-Yaacov-Berenstein-Henson-Usvyatsov:2008} suffices for our purposes.
\end{remarks}

\subsection{Tsirelson Structures}
\label{sec:Tsirelson-structures}

The \emph{classical Tsirelson structure~$\sT$} is the $\LNseq$ structure obtained from Tsirelson's space~$\cT$, with $\nrmB_k$ interpreted as the $k$-th Tsirelson approximant, and $\NrmB$ interpreted as the Tsirelson norm (via the $\NRM$ predicate, of course).

A \emph{Tsirelson structure} is a model~$\cM$ of~$\ThTs$, the $\LNseq$-theory of the classical Tsirelson structure~$\sT$.
In general, the $\bN$-sort $\bN^{\cM}$ of a Tsirelson structure~$\cM$ contains elements not named by the symbols $\NrmB$, $\nrmB_k$ (see Remark~\ref{rem:master-norm}).
However, any element $N\in\bN^{\cM}$ whatsoever still serves as the name of a seminorm $x\mapsto N(x) \coloneq \NRM(N,x)$ on the vector sort~$V^{\cM}$.
$\ThTs$ also ensures that $N(\cdot)\le\Nrm{\cdot}^{\cM}$ everywhere.

  In every Tsirelson structure $\cM$ we define the (external) \emph{Tsirelson norm~$\nrmT{\cdot}^{\cM}$} on~$\cM$ in the obvious way, namely
  \begin{equation}\label{eq:Tsirelson-external}
    \nrmT{x}^{\cM} \coloneq  \lim_{k\to\infty}\nrm{x}_k^{\cM}
    \quad\text{for all $x\in V^{\cM}$.}
  \end{equation}
  The existence (and finiteness!) of the limit in~\eqref{eq:Tsirelson-external} follows from the observation that~$\cM$ is a model of~$\ThTs$, which includes the sentences $(\forall x\in V_1)\bigl(\nrm{x}_k\dotle\nrm{x}_l\dotle\Nrm{x}\bigr)$ for $k<l\in\NN$, since these are true in the classical Tsirelson space~$\sT$ where $\NrmB^{\sT}$ coincides with~$\nrmT{\cdot}^{\sT}$.
  Thus, for all $x\in V_1^{\cM}$ (and, by linearity, for all~$x\in V^{\cM}$),
  \begin{equation*}
    \nrmT{x}^{\cM}
    = \lim_{k\to\infty}\nrm{x}_k^{\cM}
    = \sup_{k\in\NN} \nrm{x}_k^{\cM}
    \le \Nrm{x}^{\cM} < \infty.
  \end{equation*}
  However, there is no \emph{a priori} reason for $\nrmT{\cdot}^{\cM}$ to coincide with $\Nrm{\cdot}^{\cM}$ in arbitrary Tsirelson structures (see Remark~\ref{rem:nonclassical-Tsirelson} below).

  From its construction, one may regard $\nrmT{\cdot}^{\cM}$ as a real-valued predicate on~$V^{\cM}$.
  On the other hand, letting the Tsirelson norm~$\nrmT{\cdot}^{\cM}$ may also be regarded as a $\NRM$-type (i.e., a type for the formula $\NRM(N,x)$).
  In fact, $\nrmT{\cdot}^{\cM}$ is a \emph{semidefinable limit} of the types $\NRM_{\nrmB_k^{\cM}}$ in the sense that it is an accumulation point (equivalently, an ultralimit%
  \footnote{In fact, $\nrmT{\cdot}^{\cM}$ is the \emph{unique} such limit point (ultralimit) when $\cM$ models $\ThTs$.})
  of the set of these types.
  
\begin{remarks}\label{rem:nonclassical-Tsirelson}
  A model~$\cM$ of~$\ThTs$ is elementarily equivalent to the classical Tsirelson space~$\sT$.
  However, we emphasize that $\nrmT{\cdot}^{\cM}$ is not part of the structure~$\cM$, but only defined externally:
  The non-definability of the Tsirelson norm (as captured in Theorem~\ref{thm:nondef-Tsirelson} below) implies that the interpretation $\Nrm{\cdot}^{\cM}$ in such a model~$\cM$ need not coincide with the (external) pointwise limit $\nrmT{\cdot}^{\cM} \coloneq  \lim_{k\to\infty}\nrm{\cdot}^{\cM}_k$.

    To conform to the multisorted framework, the formula~$\varphi$ underlying the $\NRM$-types will be, strictly speaking, the formula $\NRM_1(N,x)$ (with $N$ of kind $\bN$ and $x$ of kind $V_1$) when issues of boundedness play a role in subsequent discussions.
\end{remarks}

\subsection{Non-Definability of the Tsirelson Norm}
\label{sec:nondef-Tsirelson}

\begin{theorem}
  \label{thm:unstab-Tsirelson}
  The Tsirelson norm is not stable in the classical Tsirelson space~$\sT$.
  More precisely, there exist:
  \begin{itemize}
  \item a sequence $\{x_i : i\in\NN\}$ in Tsirelson's space, and 
  \item a sequence $\{\nrm{\cdot}_{k_j} : j\in\NN\}$ of approximants to the Tsirelson norm~$\nrmT{\cdot}$,
  \end{itemize}
  such that
  \begin{equation}\label{eq:unequal-limits}
    \lim_{i\to\infty}\nrmT{x_i} =
    \lim_{i\to\infty}\lim_{j\to\infty}\nrm{x_i}_{k_j}
    \ne
    \lim_{j\to\infty}\lim_{i\to\infty}\nrm{x_i}_{k_j},
  \end{equation}
  where all limits involved exist and are finite.
\end{theorem}
\begin{proof}
  With $\{x_i\}$ and $\{k_j\}$ chosen so the conclusions of Proposition~\ref{prop:double-limit-Tsirelson} hold, all limits in the statement of Theorem~\ref{thm:unstab-Tsirelson} exist, and
  \begin{equation*}
   1 = \lim_{i\to\infty}\nrmT{x_i}
      = \lim_{i\to\infty}\lim_{j\to\infty}\nrm{x_i}_{k_j},
  \end{equation*}
  but
  \begin{equation*}
    \lim_{j\to\infty}\lim_{i\to\infty}\nrm{x_i}_{k_j}
    \le 1/3.\qedhere
  \end{equation*}
\end{proof}


\begin{proposition}\label{prop:definable-predicate-on-types}
Given a structure~$\cM$, the following properties are equivalent:

\begin{enumerate}
\item A formula~$\varphi$ is stable over a structure~$\cM$.
\item If a $\varphi$-type $p\in\cS_{\varphi}(\cM)$ is an accumulation point of a sequence~$(p_i)_{i=1}^{\infty}$ of $\varphi$-types realized in~$\cM$, then $p$ is definable over~$\cM$, and there exists a subsequence $(p_{i_j})_{j=1}^{\infty}$ whose $\varphi$-types converge point-wise on~$\mathcal{S}_{\tphi}(\cM)$ to a definition of~$p$.
\end{enumerate}
Furthermore, if such is the case, then every $\varphi$-type is a pointwise limit (i.e., a limit in the logic topology) of realized types: $\cS^{\cM}_{\varphi} = \overline{\mathcal{P}_{\varphi}^{\cM}}.$
\end{proposition}
\begin{proof}
   The assertion is a particular case of the equivalence between properties (i) and~(iii) in~\cite[Theorem~5]{BenYaacov2014}.
\end{proof}

\begin{theorem}[First-Order Non-Definability of the Tsirelson Norm]
  \label{thm:nondef-Tsirelson}
  The Tsirelson norm $\nrmT{\cdot}^{\cM}$ is not definable over $\LNseq$-structures~$\cM$ that satisfy either of the following properties:
  \begin{itemize}
  \item $\cM$ is an elementary extension of the classical Tsirelson space~$\sT$, or
  \item $\cM$ is a countably saturated model of~$\ThTs$.
  \end{itemize}
\end{theorem}
\begin{proof}
  As shown in section~\ref{sec:Tsirelson-structures}, $\ThTs$ implies that the Tsirelson norm~\eqref{eq:Tsirelson-external} is the only global type semidefinable over the sequence of types $\NRM_k\coloneq\nrm{\cdot}_k$ of the approximant norms (independently of the choice of~$\cU\in\beta\NN$).%
  \footnote{Recall that the relevant types and predicates are relative to the formula $\NRM_1(N,x)$ with $x$ of sort~$V_1$.
  The present use of the subindex $k$ in $\NRM_k$ should cause no confusion.}
  In particular, the Tsirelson norm is the only accumulation point of any sequence~$(\nrm{\cdot}_{k_i})_{i=1}^{\infty}$ of the $\NRM$-types of the approximant norms.
  
  Assume first that $\cM$ is an elementary extension of~$\sT$.
  By Theorem~\ref{thm:unstab-Tsirelson}, $\NRM$ is not stable in~$\sT$, hence neither in~$\cM\succeq\sT$;
  thus, by Proposition~\ref{prop:definable-predicate-on-types}, the Tsirelson norm is not definable over~$\cM$, but we can be more specific.
  With $(x_i)$ and $(k_j)$ as in Theorem~\ref{thm:unstab-Tsirelson}, let $\cU$ be any nonprincipal ultrafilter on~$\NN$, and let $q = \lim_{i,\cU}\NRM^{\sim}_{x_i}$ be the $\cU$-ultralimit type of the realized  dual types of~$(x_i)$.
  If $\nrmT{\cdot}^{\cM}$ were definable, then it would be equal to the pointwise limit of the types $\nrm{\cdot}_{l_j}$ for some subsequence $(l_j)$ of~$(k_j)$, and the value (of the extension to dual types%
  \footnote{See Remark~\ref{rem:semidef-predic-type}})
  of $\nrmT{\cdot}$ at~$q$ would be~$\lambda \coloneq  \lim_{j\to\infty}\xi_{l_j}(q) = \lim_{j\to\infty}\lim_{i,\cU}\nrm{x_i}_{l_j} = \lim_{j\to\infty}\lim_{i\to\infty}\nrm{x_i}_{k_j}$ (since the latter limits exist and equal the respective subsequential limit and ultralimit).
  On the other hand, again by definability, $\nrmT{\cdot}$ would be a continuous function on dual types, so its value at~$q = \lim_{i,\cU}\NRM^{\sim}_{x_i}$ would equal $\rho \coloneq  \lim_{i,\cU}\nrmT{q}^{\cM} = \lim_{i,\cU}\nrmT{x_i}^{\cM} = \lim_{i\to\infty}\lim_{k\to\infty}\nrm{x_i}_k^{\cM} = \lim_{i\to\infty}\lim_{j\to\infty}\nrm{x_i}_{k_j}^{\cM}$.
  However, $\lambda\ne\rho$ per Theorem~\ref{thm:unstab-Tsirelson}.
  Thus, in the sense explained, the dual type~$q$ witnesses the non-definability of~$\nrmT{\cdot}^{\cM}$.
  
  Next, let~$\cM$ be a countably saturated model of~$\ThTs$.
  Consider the full dual type~$q$ (in a countable tuple $\overline{z} = (z_i)_{i=1}^{\infty}$ of variables) of the sequence~$\overline{x} = (x_i)_{i=1}^{\infty}$ in~$V_1^{\sT}$ above.
  This full dual type is the collection $(\tphi_{\overline{x}} : \varphi\in\Phi)$ of all dual types $\tphi_{\overline{x}}  = \varphi(\cdot;\overline{x})$ as $\varphi(N;\overline{z})$ varies over the collection~$\Phi$ of all $L$-formulas in which only finitely many of the variables~$z_i$ appear free.
  The dual types~$\varphi_{\overline{x}}$ considered here are functions on the collection~$\bigl\{\Nrm{\cdot}^{\sT}\bigr\} \cup \bigl\{\nrm{\cdot}_k^{\sT}\bigr\}_{k=1}^{\infty}$ of all interpretations of the norm symbol~$N$ in~$\sT$ (rather than functions defined only on the collection of norms~$\nrm{\cdot}_k^{\sT}$).%
  \footnote{\label{fn:Nrm-not-Tsir}In $\sT$, the interpretation $\Nrm{\cdot}^{\sT}$ of the master norm is the Tsirelson norm~$\nrmT{\cdot}^{\sT}$, but this need not hold in the model~$\cM$.
    (See Remark~\ref{rem:nonclassical-Tsirelson}).}
  Since $q$ is a type on countably many variables, countable saturation of~$\cM$ ensures that it is realized by a sequence~$\overline{y} = (y_i)_{i=1}^{\infty}$ in~$V^{\cM}$.
  For ease of exposition, let us assume that~$(x_i)$ and $(k_j)$ were chosen with the properties stated in Proposition~\ref{prop:triangle-Tsirelson}.
  For each~$i$, the element~$y_i$ satisfies%
  \footnote{However, the external Tsirelson norm~$\Nrm{y_i}_T^{\cM}$ need not equal~$\Nrm{x_i}^{\sT} = \nrmT{x_i}^{\sT} = 1$.}
  \begin{itemize}
  \item $\nrm{y_i}_{k_j} \le 1/3$, if $j\le i$,
  \item $\nrm{y_i}_{k_j} \ge 2/3$, if $j > i$, and
  \item $\Nrm{y_i}^{\cM} = 1$,
  \end{itemize}
  since the corresponding properties of~$x_i$ are translated to an equality (or inequality) satisfied by the dual type~$\NRM^{\sim}_{x_i} \in q$ at the points~$N = \Nrm{\cdot}^{\cM}$ (or~$N=\nrm{\cdot}_{k_j}^{\cM}$), and $y_i$ realizes~$q$.
  $\ThTs$ ensures that, in models~$\cM$ thereof, the interpretations of the norms satisfy~$\nrm{\cdot}_k^{\cM}\le\nrm{\cdot}_l^{\cM}\le\Nrm{\cdot}^{\cM}$ pointwise for $k\le l$.
  In particular, $\nrmT{\cdot}^{\cM}$ is (finite) and everywhere defined on~$V^{\cM}$.
  In fact
  \begin{equation*}
    \nrmT{y_i}^{\cM}
    = \lim_{j\to\infty}\nrm{y_i}_{k_j}^{\cM}
    = \sup_k\nrm{y_i}_k^{\cM}
    \le \Nrm{y_i}^{\cM} = 1.
  \end{equation*}
  On the one hand, for all~$i\in\NN$, we have
  \begin{equation*}
    \nrmT{y_i}^{\cM}
    \ge \nrm{y_i}_{k_{i+1}}^{\cM}
    \ge 2/3,
  \end{equation*}
  hence
  \begin{equation*}
    \lim_{i\to\infty}\lim_{j\to\infty}\nrm{y_i}_{k_j}
    = \lim_{i\to\infty}\nrmT{y_i}^{\cM}
    \ge 2/3.
  \end{equation*}
  Furthermore, for fixed~$j\in\NN$, the sequence $\{\nrm{y_i}_{k_j}\}_{i=1}^{\infty}$ takes values $\le 1/3$ for $i\ge j$;
  thus, $\lim_{i\to\infty}\nrm{y_i}_{k_j} \le 1/3$, and hence
  \begin{equation*}
    \lim_{j\to\infty}\lim_{i\to\infty}\nrm{y_i}_{k_j}
    \le 1/3.
  \end{equation*}
  This shows that the formula~$\NRM$ is not stable in~$\cM$, and we conclude that the Tsirelson norm is not definable over~$\cM$, by Theorem~\ref{thm:unstab-Tsirelson}.
\end{proof}

\section{Schlumprecht's space}
\label{sec:Schlumprecht}

\subsection{Construction of the Schlumprecht space}
\label{sec:Schlumprecht-cons}

The construction of Schlumprecht's space is similar to that of Tsirelson's.
We use the same notation as in section~\ref{sec:Tsirelson}, except for the requirement that admissible families of finitely many nonempty finite subsets, say $E_i$ ($1\le i\le m$) of~$\NN$ need only satisfy $E_1<E_2<\dots<E_m$ (and not necessarily the requirement $E_1\ge m$ as in Tsirelson's construction).
Let $L(t) = \log_2(t+1)$ and $\fL(t) = t/L(t)$ for $t>0$.
It is trivial to verify that both $L$ and $\fL$ are unbounded and strictly increasing.

Define the sequence $\{\nrm{\cdot}_k\}_{k=0}^{\infty}$ of norms on~$c_{00}$ (the \emph{Schlumprecht approximants}) recursively as follows.
For $x = \sum_na_ne_n\in c_{00}$,
\[
  \nrm{x}_0 \coloneq  \nrminf{x} = \max_n|a_n|,
\]
and, for $k\ge 1$,
\[
  \nrm{x}_k =
  \max\left\{
    \frac{1}{L(m)}\sum_{i=1}^m\nrm{E_ix}_{k-1}:
    E_1 < E_2 <\dots < E_m
  \right\},
\]
where the maximum is taken over the set of all admissible collections of any size $m\ge 1$.
(The reuse of the notation $\nrm{\cdot}_k$ is convenient, but unrelated to the Tsirelson norms introduced in section~\ref{sec:Tsirelson-defns}.)
Since $L(1) = 1$, we see that the inequality $\nrm{x}_k \ge \nrm{x}_{k-1}$ holds for all $x\in c_{00}$ and $k\ge 1$, as witnessed by the singleton collection $\{E_1\}$, where $E_1=\supp(x)$.
The \emph{Schlumprecht norm} of $x\in c_{00}$ is defined by
\[
  \nrm{x} = \lim_{k\to\infty}\nrm{x}_k\quad(=\sup_{k\in\NN}\nrm{x}_k),
\]
and \emph{Schlumprecht's space~$\cS$} is the norm-completion of~$(c_{00},\nrmS{\cdot})$.
The induced norm on $\cS$ is still denoted $\nrmS{\cdot}$ by an abuse of notation.
Routine induction shows that $\nrm{x}_k\le \sum_n|a_n| = \nrmone{x}$ for $x\in c_{00}$ and $k\in\NN$.
Thus, inequality $\nrmS{x} \le \nrmone{x}$ also holds for all $x\in c_{00}$ \emph{a priori} and, \emph{a posteriori}, for all $x\in\cS$.

Evidently, the Schlumprecht approximants and the Schlumprecht norms depend only on the absolute values of the coefficients, and monotonically so.

In contrast to Tsirelson's construction, the admissibility condition for families $\{E_i\}$ in Schlumprecht's case is invariant under the shift-by-$j$ transformation $\sum_ia_ie_i \mapsto \sum_ia_ie_{i+j}$ for any fixed $j\ge 0$.
Since $\nrm{\cdot}_0 = \nrminf{\cdot}$ is also shift-invariant, it follows inductively that all approximants $\nrm{\cdot}_k$ are shift-invariant, and so is the Schlumprecht norm $\nrmS{\cdot}$ itself.%
\footnote{More generally, the Schlumprecht norm is invariant under monotone re-indexing transformations $\sum_i a_i e_i \mapsto \sum_i a_i e_{j_i}$ via any fixed choice of strictly increasing indexes $1\le j_1 < j_2 < \dots < j_i < \dots$.}

\subsection{Non-uniform convergence of the Schlumprecht approximants.}
\label{sec:nonunif-Schlum}

\subsubsection{Auxiliary sequences of rapid growth}
\label{sec:Schl-seqs}

\begin{definition}\label{def:aux-fns-Schlum}
  Let $L^{*}$, $\fL^{*}$ be the integer-valued quasi-inverses of $L$, $\fL$, namely, for $n\in\NN$,
  \begin{align*}
    \fL^{*}(n)
    &= \lceil \fL^{-1}(n)\rceil
      = \min\{m\in\NN : m/\log_2(m+1) \ge n\}; \\
    L^{*}(n)
    &= L^{-1}(n)
      = 2^n-1
      \quad (= \min\{m\in\NN : \log_2(m+1)\ge n\}).
  \end{align*}
\end{definition}

\begin{proposition}\label{prop:aux-fns-Schl}
  Given any function $M:\NN\to\NN$, there exist unique functions
  \begin{align*}
    \Lambda &: (k,n) \mapsto \Lambda_k(n)\quad(k,n\ge 1) \\
    \lambda &: (k,n,i) \mapsto \lambda_k^i(n)\quad(k,n,i\ge 1) \\
    \nu &: (k,n,i) \mapsto \nu_k^i(n)\quad(k\ge 2\ \&\ n,i\ge 1)
  \end{align*}
  taking values in~$\NN$ such that the following identities hold:
  \begin{enumerate}
  \item $\lambda_1^1(n) = 1$;
  \item $\lambda_1^{i+1}(n) = L^{*}\bigl(2L[\lambda_1^i(n)]\bigr) = [1+\lambda_1^i(n)]^2-1$;
  \item $\lambda_k^i(n) = \Lambda_{k-1}[2^i\nu_k^i(n)]$\quad($k\ge 2$);
  \item $\Lambda_k(n) = \sum_{i=1}^{M(n)}\!\lambda_k^i(n)$;
  \item $\nu_k^1(n) = 2$\quad ($k\ge 2$);
  \item $\nu_k^{i+1}(n) = L^{*}\bigl(2^iL[\lambda_k^i(n)]\bigr) = [1+\lambda_k^i(n)]^{2^i}-1$\quad($k\ge 2$).
\end{enumerate}
If $M$ is strictly increasing 
then each of the functions $\lambda,\Lambda, \mu$ is increasing in each variable separately (strictly increasing indeed---apart from the equalities $\lambda_1^1(n)= 1$ and $\nu_k^1(n) = 2$ for all $k\ge 2$ and $n\ge 1$).
\end{proposition}
\begin{proof}
  The existence and uniqueness of $\Lambda$, $\lambda$ and~$\nu$ are immediate, since the given conditions amount to a jointly recursive definition thereof.
  The asserted monotonicity of $\lambda,\Lambda,\mu$ conditional on $M$ being strictly increasing is trivially verified.
\end{proof}

 For $k,n,i\ge 1$, define
  \begin{equation*}
    \begin{split}
      \lambda_k^{\les 0}(n) &= 0,\qquad\text{and}\\
      \lambda_k^{\les i}(n) &= \sum_{j=1}^i\lambda_k^j(n).
    \end{split}
  \end{equation*}
  (In particular, $\Lambda_k(n) = \lambda_k^{\les M(n)}(n)$ for $k\ge 1$.)


  \subsubsection{The main proposition}
  \label{sec:Schlum-propn}

  \begin{proposition}
    \label{prop:Schlumprecht-bounds}  
    There exist unique vectors $Z(l)$, $X_k(n)$, $Y_k(n)$ and $x_k^i(n)$ ($k,l,n,i\ge 1$) in $c_{00}$ such that
    \begin{enumerate}
    \item $Z(l) = \frac{1}{\fL(l)}\sum_{i=1}^le_i$;
    \item $X_k(n) = Y_k(n)/\nrm{Y_k(n)}_{k+1}$;
    \item $Y_k(n) = \mathrm{D}\sum_{i=1}^{M(n)}x_k^i(n)$;
    \item $x_1^i(n) = Z\bigl(\lambda_1^i(n)\bigr)$;\quad and
    \item $x_k^i(n) = X_{k-1}\bigl(2^i\nu^i_k(n)\bigr)$\quad ($k\ge 2$);
    \item The support of $Y_k(n)$ is $I_k(n) \coloneq  \bigl[1,\Lambda_k(n)\bigr]$;
    \item The support of $x_k^i(n)$ is $\bigl[1, \lambda_k^i(n) \bigr]$, and the coefficients of $x_k^i(n)$ in the direct sum~(3) fill the interval $I_k^i(n) \coloneq  \bigl(\lambda_k^{\les i-1}(n), \lambda_k^{\les i}(n) \bigr]$ in $Y_k(n)$, and also in~$X_k(n)$;
    \item $\nrm{x_k^i(n)}_k = 1$;
    \item $\nrm{Y_k(n)}_k \le 3$;\quad and
    \item $\nrm{Y_k(n)}_{k+1} \ge 3n$.
    \end{enumerate}
  \end{proposition}

  To begin the proof of Proposition~\ref{prop:Schlumprecht-bounds}, note first that conditions~(1)--(5) amount to a (unique) recursive definition of all the required vectors.
    
    Properties~(6) and~(7) are proved simultaneously by induction on $k$ (and, for $k$ fixed, by induction on~$i$).
    The details are trivial and omitted.
    
    We prove that~(8), (9) and~(10) hold, for all $n\ge 2$ and $i\ge 1$, by induction on~$k\ge 1$.
Throughout the proof, we write $M$ for $M(n)$.
    
\subsubsection{The base of the induction}

The proof for $k=1$ is as follows.
We will presently write $Y(n)$, $I_j$, $x_j$, $\lambda_j$ for $Y_1(n)$, $I_1^j(n)$, $x_1^j(n)$, $\lambda_1^j(n)$ (the last three expressions are constant in~$n$).
We remark that, in the definition of~$\nrm{\sum_{i=1}^le_i}_1$, the admissible family achieving the maximum value, equal to $l/L(l) = \fL(l)$, is $E_i=\{i\}$ ($1\le i\le L$);
therefore, $\nrm{Z(l)}_1 = 1 = \nrm{x^i}_1$ for $l,i\ge 1$ (by~(1) and~(4)), proving~(8) for $k=1$.

Using the family $\{I^i\}_{i=1}^{M}$ in the definition of $\nrm{\cdot}_2$, we have $\nrm{E_iY(n)}_1 = \nrm{x^i}_1=1$ (by translation invariance of $\nrm{\cdot}_1$), hence the lower bound
\begin{equation*}
  \nrm{Y(n)}_2
  \ge \frac{1}{L(M)}\sum_1^{M}1
  = \frac{M}{L(M)} = \fL(M) \ge 3n,
\end{equation*}
since $M = \fL^{*}(3n)$.
This proves~(10).

To find an upper bound for $\nrm{Y(n)}_1$, note first that the coefficients of the basis elements $e_j$ in $Y(n)$ decrease with~$j$ (since the sequence $1/\fL\bigl(\lambda^i\bigr)$ decreases with~$i$).
By the definitions of $\nrm{\cdot}_0$ and $\nrm{\cdot}_1$, a moment's reflection shows that a family realizing the maximum that defines $\nrm{Y(n)}_1$ must consist of, say, $l$~singletons $E_i = \{i\}$ for some $l\le\Lambda(n)$ and $1\le i\le l$.
Choose $m \le M$ so that $l\in(\lambda^{\les m-1}, \lambda^{\les m}]$, and let $l' = l-(\lambda^{\les m-1}) \le \lambda^m$.
Then, we have
\begin{equation}
  \label{eq:nrm1-Y1}
  \begin{split}
    \nrm{Y(n)}_1
    &= \frac{1}{L(l)}
    \left(
      \sum_{i=1}^{m-1}\sum_{j=1}^{\lambda^i}x^i[j]
      + \sum_{j=1}^{l'}x^m[j]
    \right) 
    = \frac{1}{L(l)}
    \left(
      \sum_{i=1}^{m-1}
      \frac{\lambda^i}{\fL(\lambda^i)}
      +
      \frac{l'}{\fL(\lambda^m)}
    \right)\\
    &=
      \sum_{i=1}^{m-1}\frac{L(\lambda^i)}{L(l)}
      +
      \frac{l'}{L(l)\fL(\lambda^m)}
    \le \sum_{i=1}^{m-1}\frac{L(\lambda^i)}{L(\lambda^{m-1})}
    +
    \frac{\fL(l')}{\fL(\lambda^m)}\qquad\text{(since $\lambda^i,l'\le l$)}\\
    &\le \sum_{i=1}^{m-1}\frac{1}{2^{m-i-1}}
    + 1
    \le 2+1=3\qquad\text{(since $l'\le \lambda^m$ and $L(\lambda^{i+1})\ge 2L(\lambda^i)$)}.
  \end{split}
\end{equation}
This proves~(9) and completes the proof of the base case $k=1$ of the induction.

\subsubsection{The inductive step}

Now we carry out the inductive step of the proof of~(8), (9) and~(10).
Assume they hold for some $k\ge 1$.
Denote $\lambda_{k+1}^i(n)$, $\nu_{k+1}^i(n)$, $x_{k+1}^i(n)$, $X_{k+1}(n)$, $Y_{k+1}(n)$ and $I^i_{k+1}(n)$ by $x_i$, $X$, $Y$, $I_i$, $\lambda_i$ and~$\nu_i$.
First, observe that $\nrm{X_k(\cdot)}_{k+1} = 1$ follows from the inductive hypothesis~(2);
therefore, $\nrm{x_i}_{k+1} = 1$ follows from~(5), proving~(8) for $k+1$.

The intervals $I_i$ ($1\le i\le M$) satisfy $\nrm{I_iY}_{k+2} = \nrm{x_i}_{k+2} = 1$.
Using these intervals in the definition of $\nrm{\cdot}_{k+2}$ we obtain
\begin{equation*}
  \nrm{Y_{k+1}(n)}_{k+2} 
 \ge
  \frac{1}{L(M)} \sum_{i=1}^{M}\nrm{x^i_{k+1}(n)}_{k+1}
  = \frac{M}{L(M)}
  = \fL(M) \ge 3n.
\end{equation*}
This proves~(10) for $k+1$.
 
To prove~(9) for $k+1$, we start as in the proof of the base case.
Since the coefficients of~$Y$ are decreasing, a moment's reflection shows that, an admissible family achieving the maximum that defines $\nrm{Y_{k+1}(n)}_{k+1}$ may be taken to consist of \textit{intervals} $E_1<E_2<\dots<E_m$.
By the monotonicity of $\nrm{\cdot}_k$, the maximizing $E_j$'s may be taken to be the back-to-back intervals such that $\bigcup_{i=1}^mE_i = [1,\Lambda] = \bigcup_{j=1}^M I_j$.

A pair $(j,i)$ with $1\le j\le m$ and $1\le i\le M$ will be called \emph{relevant} if $E_j\cap I_i$ is nonempty.

We split the intervals $[1,m]$ into two disjoint subclasses:
\begin{itemize}
\item[($\subseteq$)]
  This class consists of those $j$ such that $E_j\subseteq I_i$ for some~$i$.
\item[($\not\subseteq$)]
  This is the complementary class consisting of those relevant $j$ such that for no~$i$ is $E_j$ a subset of $I_i$.
\end{itemize}
We will use the symbols ``$\subseteq$'', ``$\not\subseteq$'' as nomenclature for the corresponding class.
  Abusing the nomenclature, we will say that $i\in[1,M]$ \emph{is of class~$\subseteq$} if $(j,i)$ is of class~$\subseteq$ for some~$j$.
With these notations, we have
\begin{equation}
  \label{eq:sum-three-cases}
\nrm{Y}_{k+1} = \sumsube + \sumnotsube,
\end{equation}
  where
 \begin{align*}
   \sumsube
   &\coloneq \sumsube_j \frac{\nrm{E_jY}_k}{L(m)}, &
\sumnotsube
   &\coloneq \sumnotsube_j \frac{\nrm{E_jY}_k}{L(m)}.
  \end{align*}
We will find an upper bound for each of the two sums above.

\emph{Case~$\subseteq$:}\quad
Consider any fixed $i\le M$ of class~$\subseteq$.
Let $m_i$ be the number of indexes~$j$ of class~$\subseteq$ for~$i$;
they form an interval, say $[r_i, r_i+m_i-1]$.
Since such $E_j$ are disjoint nonempty subsets of~$I_i$, and $\#I_j = \lambda_j$, we have $m_i\le\min\{M,\lambda_i\}$.
Now, we have
  \begin{equation}
    \label{eq:subseteq-reverse-order}
    \sumsube
    =
    \sumsube_j \frac{\nrm{E_jY}_k}{L(m)}
    =
    \sumsube_i
    \sum_{j=r_i}^{r_i+m_i-1}
    \frac{\nrm{E_jY}_k}{L(m)}.
  \end{equation}
as seen by using the family $\{E_j\}_{j=s_i}^{s_i+m_i-1}$ in the recursive definition of $\nrm{I_1Y}_{k+1}$, which is admissible precisely because such $(i,j)$ are of type~$\subseteq$.

If there exists $i\in[1,M]$ of class~$\subseteq$ such that $m_i > \nu_i$, let $p$ denote the largest such~$j$;
otherwise, let $p=0$.
  Thus, every $i>p$ of class~$\subseteq$ satisfies $m_j \le \nu_j$.
Continuing from~\eqref{eq:subseteq-reverse-order}, we write
\begin{equation}
  \label{eq:short-long}
  \sumsube = \sumsube_{\les p} + \sumsube_{>p},
\end{equation}
where $\sumsube_{\les p}$, $\sumsube_{>p}$ are the sums over (class-$\subseteq$) indices $i\le p$, $i>p$, respectively, on the right-hand side of~\eqref{eq:subseteq-reverse-order}.

If $p>0$, we have
\begin{equation}\label{eq:long-estimate}
  \begin{split}
    \sumsube_{\les p}
    &\coloneq  \sumsube_{i\le p}\sum_{j=r_i}^{r_i+m_i-1} \frac{\nrm{E_jY}_k}{L(m)}
    \le \sumsube_{i\le p}
    \frac{L(m_i)\nrm{I_iY}_{k+1}}{L(m)}\\
    &\qquad\text{(by definition of~$\nrm{\cdot}_{k+1}$, since $\{E_j\}_{j=r_i}^{r_i+m_i-1}$ is an admissible family of subsets of~$I_i$)}\\
    &= \sumsube_{i\le p}
    \frac{L(m_i)\nrm{x_i}_{k+1}}{L(m)}
    = \sumsube_{i\le p}
    \frac{L(m_i)}{L(m)}
    \qquad\text{(by property~(8)---already proved for $k+1$)}\\
    &\le 1 + \sum_{i=1}^{p-1} \frac{L(\lambda_i)}{L(\nu_p)}
    \qquad\text{(since $m_i\le\lambda_i$ and $m\ge m_p > \nu_p$, by choice of~$p$)}\\
    &\le 1 + \sum_{i=1}^{p-1} \frac{L(\lambda_i)}{L(\nu_{i+1})}
    = 1 + \sum_{i=1}^{p-1} \frac{L(\lambda_i)}{L(\nu_{i+1})}
    = 1 + \sum_{i=1}^{p-1}\frac{1}{2^i}\\
    &\qquad\text{(since $\nu_i < \nu_p$ for $i<p$, and $L(\nu_{i+1}) = 2^iL(\lambda_i)$ by construction).}
  \end{split}
\end{equation}
We have $\sumsube_{\les p} = 0$ if $p=0$, i.e., if there are no class-$\subseteq$ indexes~$i$ such that $m_i>\nu_i$, so the upper bound (equal to~$1$) on the right-hand side of inequality~\eqref{eq:long-estimate} remains valid in this case as well.

We now estimate $\sumsube_{>p}$.
Just as in the case of $\sumsube_{\le p}$ above, have $\sumsube_{>p} = 0$ if there is no index $i>p$ (perhaps no index~$i$ at all) of class~$\subseteq$.
  Letting $n_i = 2^i\nu_i$, we have
\begin{equation}
\label{eq:nrmk-x-k+1}
\nrm{x_i}_k 
= \nrm{X_k(n_i)}_k
= \frac{\nrm{Y_k(n_i)}_k}
{\nrm{Y_k(n_i)}_{k+1}} \le \frac{3}{3n_i}
= \frac{1}{2^i\nu_i}.
\end{equation}
We get,
\begin{equation}
  \label{eq:short-estimate}
  \begin{split}
        \sumsube_{>p}
        &\coloneq  \sumsube_{i>p}\sum_{j=r_i}^{r_i+m_i-1} \frac{\nrm{E_jY}_k}{L(m)}
        = \sumsube_{i>p}\sum_{j=r_i}^{r_i+m_i-1} \frac{\nrm{I_iY}_k}{L(m)}
        = \sumsube_{i>p}\frac{m_i\nrm{x_i}_k}{L(m)}\\
        &\qquad\text{(by monotonicity of $\nrm{\cdot}_k$, since $E_j\subseteq I_i$; furthermore, $I_iY$ is a shift of $x_i$)}\\
        &\le
        \sumsube_{i>p}\frac{\nu_i\nrm{x_i}_k}{L(m)}
        \le
        \sumsube_{i>p}\frac{2^{-i}}{L(m)}
        \le \sum_{i=p+1}^{\infty}\frac{1}{2^i}
        \qquad\text{(by~\eqref{eq:nrmk-x-k+1}, since $m_i\le \nu_i$ for $i>p$ of class~$\subseteq$).}
  \end{split}
\end{equation}

Combining estimates~\eqref{eq:long-estimate} and~\eqref{eq:short-estimate} we obtain
\begin{equation}
  \label{eq:bound-sumsube}
  \sumsube \le 2.
\end{equation}

\emph{Case~$\not\subseteq$}:\quad
For $j$ of class $\not\subseteq$, the relevant indexes $i$ fill an interval $J_j = [s_j, s_j+n_j-1]$ of length $n_j\ge 2$.
For each $i\in[1,M]$, let $\mu_i$ be the number of $j$ of class~$\not\subseteq$ such that $(j,i)$ is relevant.
Clearly, $0\le \mu_i\le 2$.
Define $F_{ji} = E_j\cap I_i$.
By the triangle inequality and monotonicity,
\begin{equation}
  \label{eq:triangle-notsube}
  \begin{split}
    \sumnotsube
    &\coloneq  \sumnotsube_j \frac{\nrm{E_jY}_k}{L(m)}
    \le \sumnotsube_j\sum_{i=s_j}^{s_j+n_j-1} \frac{\nrm{F_{ji}Y}_k}{L(m)}
    \le \sumnotsube_j\sum_{i=s_j}^{s_j+n_j-1} \frac{\nrm{I_iY}_k}{L(m)}\\
    &= \sum_i\mu_i \frac{\nrm{I_iY}_k}{L(m)}
    = 2\sum_i \frac{\nrm{x_i}_k}{L(m)}
    \le \frac{2}{L(m)}\sum_{i=1}^M \frac{1}{2^i\nu_i}
    \le \frac{2}{\nu_1} = 1.
  \end{split}
\end{equation}
By~\eqref{eq:sum-three-cases}, \eqref{eq:bound-sumsube} and~\eqref{eq:triangle-notsube} we have proved~(9) for $k+1$, completing the inductive step and the proof of Proposition~\ref{prop:Schlumprecht-bounds}.

\subsection{Non-Definability of the Schlumprecht Norm}
\label{sec:nondef-Schlump}

\begin{theorem}
  Lemma~\ref{lem:Tsirelson-gaps}, Propositions~\ref{prop:triangle-Tsirelson} and~\ref{prop:double-limit-Tsirelson}, and Theorem~\ref{thm:unstab-Tsirelson} hold for the Schlumprecht norm~$\nrmS{\cdot}$ and its approximants $\nrm{\cdot}_k$ in place of $\nrmT{\cdot}$ and~$\nrm{\cdot}_k$ therein.
\end{theorem}
\begin{proof}
  With the indicated one-for-one replacements, the proofs go through verbatim provided at the beginning of the proof of the revised Lemma~\ref{lem:Tsirelson-gaps} we take $y = Y_k(3)/3$ as per Proposition~\ref{prop:Schlumprecht-bounds} (with $n=3$), whose Schlumprecht approximants satisfy $\nrm{y}_k\le 1$ and $\nrm{y}_{k+1}\ge 3$.
\end{proof}

The \emph{classical Schlumprecht structure~$\sS$} is the Nseq structure obtained from Schlumprecht space~$\cS$ by interpreting $\nrmB_k$ as the Schlumprecht approximants, and the master norm $\NrmB$ as the Schlumprecht norm on~$\cS$.
$\ThSch$ is the $\LNseq$-theory of~$\sS$.
A \emph{Schlumprecht structure} is any model of~$\ThSch$.

The following analogue of Theorem~\ref{thm:nondef-Tsirelson} holds.
Since the analogue of Theorem~\ref{thm:unstab-Tsirelson} (for the Schlumprecht norm) holds, the proof is \emph{mutatis mutandis} the same, and omitted.

\begin{theorem}[First-Order Non-Definability of the Schlumprecht Norm]
  \label{thm:nondef-Schlumprecht}
  The Schlumprecht norm $\nrmT{\cdot}^{\cM}$ is not definable over $\LNseq$-structures~$\cM$ that satisfy either of the following properties:
  \begin{itemize}
  \item $\cM$ is an elementary extension of the classical Schlumprecht space~$\sS$, or
  \item $\cM$ is a countably saturated model of~$\ThSch$.
  \end{itemize}
\end{theorem}


\begin{thebibliography}{BYBHU08}

\bibitem[Bel84]{Bellenot:1984}
Steven~F. Bellenot.
\newblock The {B}anach space {$T$} and the fast growing hierarchy from logic.
\newblock {\em Israel J. Math.}, 47(4):305--313, 1984.

\bibitem[BY14]{BenYaacov2014}
Ita\"{\i} Ben~Yaacov.
\newblock Model theoretic stability and definability of types, after {A}.
  {G}rothendieck.
\newblock {\em Bull. Symb. Log.}, 20(4):491--496, 2014.

\bibitem[BYBHU08]{Ben-Yaacov-Berenstein-Henson-Usvyatsov:2008}
Ita{\"{\i}} Ben~Yaacov, Alexander Berenstein, C.~Ward Henson, and Alexander
  Usvyatsov.
\newblock Model theory for metric structures.
\newblock In {\em Model theory with applications to algebra and analysis.
  {V}ol. 2}, volume 350 of {\em London Math. Soc. Lecture Note Ser.}, pages
  315--427. Cambridge Univ. Press, Cambridge, 2008.

\bibitem[BYU10]{Ben-Yaacov-Usvyatsov:2010}
Ita{\"{\i}} Ben~Yaacov and Alexander Usvyatsov.
\newblock Continuous first order logic and local stability.
\newblock {\em Trans. Amer. Math. Soc.}, 362(10):5213--5259, 2010.

\bibitem[CS89]{Casazza-Shura:1989}
Peter~G. Casazza and Thaddeus~J. Shura.
\newblock {\em Tsirel\cprime son's space}, volume 1363 of {\em Lecture Notes in
  Mathematics}.
\newblock Springer-Verlag, Berlin, 1989.
\newblock With an appendix by J. Baker, O. Slotterbeck and R. Aron.

\bibitem[FJ74]{Figiel-Johnson:1974}
Tadeusz Figiel and William~B. Johnson.
\newblock A uniformly convex {B}anach space which contains no {$l\sb{p}$}.
\newblock {\em Compositio Math.}, 29:179--190, 1974.

\bibitem[Gai76]{Gaifman:1976}
Haim Gaifman.
\newblock Models and types of {P}eano's arithmetic.
\newblock {\em Ann. Math. Logic}, 9(3):223--306, 1976.

\bibitem[Gow95]{Gowers:1995}
W.~Timothy Gowers.
\newblock Recent results in the theory of infinite-dimensional {B}anach spaces.
\newblock In {\em Proceedings of the International Congress of Mathematicians,
  Vol.\ 1, 2 (Z\"urich, 1994)}, pages 933--942, Basel, 1995. Birkh\"auser.

\bibitem[Jam64]{James1964}
R.~C. James.
\newblock Uniformly non-square {B}anach spaces.
\newblock {\em Ann. of Math.}, 80:542--550, 1964.

\bibitem[Kei]{Keisler2020}
H.~Jerome Keisler.
\newblock Model theory for real-valued structures.
\newblock To appear in {\em Beyond First-Order Model Theory, Vol.\ 2}.

\bibitem[Ode02]{Odell:2002}
Edward Odell.
\newblock On subspaces, asymptotic structure, and distortion of {B}anach
  spaces; connections with logic.
\newblock In {\em Analysis and logic (Mons, 1997)}, volume 262 of {\em London
  Math. Soc. Lecture Note Ser.}, pages 189--267. Cambridge Univ. Press,
  Cambridge, 2002.

\bibitem[Sch91]{Schlumprecht:1991}
Thomas Schlumprecht.
\newblock An arbitrarily distortable {B}anach space.
\newblock {\em Israel J. Math.}, 76(1-2):81--95, 1991.

\bibitem[She71]{Shelah:1971a}
Saharon Shelah.
\newblock Stability, the f.c.p., and superstability; model theoretic properties
  of formulas in first order theory.
\newblock {\em Ann. Math. Logic}, 3(3):271--362, 1971.

\bibitem[She90]{Shelah:1990}
Saharon Shelah.
\newblock {\em Classification theory and the number of nonisomorphic models},
  volume~92 of {\em Studies in Logic and the Foundations of Mathematics}.
\newblock North-Holland Publishing Co., Amsterdam, second edition, 1990.

\bibitem[Tsi74]{Tsirelson:1974}
Boris~S. Tsirel'son.
\newblock It is impossible to imbed $l\sb{p}$ of $c\sb{0}$ into an arbitrary
  {B}anach space.
\newblock {\em Funkcional. Anal. i Prilo\v zen.}, 8(2):57--60, 1974.

\end{thebibliography}

\def\cprime{$'$}

\end{document}